\documentclass[11pt]{article}

\usepackage{amssymb,amsthm, graphicx, stmaryrd, fancyhdr}
\usepackage{amsfonts}
\usepackage{amsmath}
\usepackage[english,german]{babel}
\usepackage{epsfig,epstopdf}
\usepackage{subfigure}
\usepackage{latexsym}
\usepackage{color}

\definecolor{darkgreen}{rgb}{0.00,0.40,0.00}

\newcommand{\ba}{\begin{array}{ll}}
\newcommand{\bal}{\begin{array}{ll}}
\newcommand{\ea}{\end{array}}

\newcommand{\E}{\mathbb{E}}

\newcommand{\Q}{\mathbb{Q}}

\newcommand{\be}{\begin{eqnarray*}}
\newcommand{\ee}{\end{eqnarray*}}
\newcommand{\ben}{\begin{eqnarray}}
\newcommand{\een}{\end{eqnarray}}

\def\sS{{\cal S}}
\def\E{\mathbb{E}}
\def\P{\mathbb{P}}
\def\Q{\mathbb{Q}}

\newtheorem{lemma}{Lemma}[section]
\newtheorem{theorem}{Theorem}[section]

\newtheorem{corollary}[theorem]{Corollary}
\newtheorem{definition}{Definition}[section]

\begin{document}

\selectlanguage{english}

\renewcommand{\baselinestretch}{1.2} \normalsize

 \title{Utility Indifference Pricing: A Time Consistent Approach \thanks{  Work supported by NSERC grant 371653-09, MITACS grant 5-26761 and the NSFC grant 10901086.}}

\author{\normalsize    Traian A.~Pirvu \\[8pt]
        \small Dept of Mathematics \& Statistics\\
        \small McMaster University \\
        \small 1280 Main Street West \\
        \small Hamilton, ON, L8S 4K1\\
        \small tpirvu@math.mcmaster.ca
        \and
        \normalsize Huayue Zhang \\[8pt]
        \small Dept of Finance\\
        \small  Nankai University\\
        \small 94 Weijin Road \\
        \small Tianjin, China, 300071 \\
        \small  hyzhang69@nankai.edu.cn
%        \small  Dept of Mathematics  \\
%        \small  Humboldt University Berlin \\
%        \small  Unter den Linden 6\\
%        \small  10099 Berlin \\
%        \small  gnreis@math.hu-berlin.de
\vspace*{0.8cm}}

\maketitle

\begin{abstract}
\thispagestyle{empty} { This paper considers the optimal portfolio
selection problem in a dynamic multi-period stochastic framework
with regime switching. The risk preferences are of  exponential (CARA) type
with an absolute coefficient of risk aversion which changes with the regime.
The market model is incomplete and there are two risky assets: one tradable and one non-tradable.
In this context, the optimal investment strategies are time inconsistent.
Consequently, the subgame perfect equilibrium strategies are considered. 
The utility indifference prices of a contingent claim written on the risky assets
are computed via an indifference valuation algorithm. By running numerical experiments,
we examine how these prices vary in response to changes in model parameters.}
\end{abstract}


\section{Introduction}
One of the most important problems in mathematical finance
is the valuation of contingent claims in incomplete financial markets.

This paper studies the indifference valuation of contingent claims in a multi-period stochastic 
model under regime switching. The risk preferences are of exponential
type and they are allowed to change with the regime.

The problem of pricing contingent claims by utility indifference
 in an incomplete binomial model was studied in \cite{MZ1} and \cite{MZ2}. The work of \cite{MZ1}
 constructs a probabilistic iterative algorithm to obtain utility indifference
 prices of contingent claims. This algorithm at each step consists of a nonlinear
 pricing functional which is applied to prices obtained at the earlier steps. This functional
 is represented in terms of risk aversion and a special martingale measure. In \cite{MZ2}, a more general
 model is considered, with an stochastic factor which may affect the transition probabilities and
 the contingent claim's payoff. Two pricing algorithms are proposed in this paper to produce the
 utility indifference prices. They employ two martingale measures: the minimal martingale measure 
 and the minimal entropy measure. This paper also analyses the dependence of the utility indifference
 prices on the choice of the trading horizon.
 
 Our paper proposes a model with regime switching. Recently, many papers considered the pricing of contingent claims
 on regime switching market. Here we recall only two such works, \cite{Guo} and \cite{Eli}. In \cite{Guo}, the author considers a stock price model which
 allows for the drift and the volatility coefficients to switch according to two-states. This market is incomplete, but it is completed with
 new securities. In \cite{Eli} the problem of option pricing is considered in a model where the risky underlying assets are driven by Markov-modulated
 Geometric Brownian motions. A regime switching Esscher transform is used to find a martingale pricing measure.
 
 The novelty our model brings is the change in risk preference during the investemnt horizon. The issue of loss aversion changing with time was addressed in financial economic literature. For instance, \cite{Bar} considers a model in which the loss aversion depends on prior gains and losses, so it
 may change through time. We choose to model this effect by allowing the risk aversion to change between two exponential type utilities according to the two
 states of the market (bull and bear). In a bull market we expect investors to be willing to take more risk and this is modeled by a lower coefficient of
 relative risk aversion as compared with the bear market. 
 
 This type of risk preferences lead to time inconsistent investemnt strategies. That is, an investor may have an incentive to deviate
 from the optimal strategies that he/she computed at some past time. To deal with this issue, \cite{Bjo} developed a theory for stochastic control problems
 which are time inconsistent in the sense that they do not admit a Bellman optimality principle. Inspired by \cite{EKE} and \cite{EKE1}, the work of \cite{Bjo} introduced the subgame perfect Nash equilibrium strategies in a discrete time model. These strategies are optimal to be implemented in the next time interval given that they are optimal
 in the future. In dealing with the problem of time consistency we choose the approach proposed by \cite{Bjo}.
 
 Our paper proposes an indifference valuation algorithm for pricing contingent claims in a discrete time incomplete market with
 regime switching. At each step, the pricing functional depends on the risk aversion, a martingale measure (the minimal martingale measure)
 and a process which keeps track of the previous optimal wealth levels. In the special case of non-switching preferences our results recover some of the results of \cite{MZ2}.
 Moreover, in this case we show that the subgame perfect equilibrium strategies coincide with the optimal ones.   
 
 The goal of this paper is to understand the impact a change in risk aversion has on prices of contingent claims. We see in times of financial crises (bear markets)
 big changes in prices and we try to explain it by an increase in investor's risk aversion. We run numerical experiments which show that an increase of
 20$\%$ in risk aversion leads to an increase in prices anywhere between $7\%$ and $23\%.$
 
 The mathematics we use to prove our results relies on a special form of optimality principle applied to a sequence of optimization problems (one for each time period).
 This involves the construction of a sequence of one period BS$\Delta$Es (one for each optimization problem). The terminal value of a BS$\Delta$E was computed in a previous time step. For the special case of non-switching preferences, these one period  BS$\Delta$Es couple in one multiperiod  BS$\Delta$E.
 In the first step we express the equilibrium strategies and utility indifference prices in terms of these  BS$\Delta$Es. In a second step, due to the special structure
 of our model, we manage to get formulas which depend only on risk aversion, the minimal martingale measure and a special process ( defined through the values of optimal wealth process at antecedent time steps).
  
  By the best of our knowledge this paper is the first work to address the indifference valuation in regime switching incomplete market with regime switching preferences.

The paper is organized as follows: Section 2 presents the model.
 Section 3 introduces the subgame perfect equilibrium strategies
 and the indifference valuation algorithm. In Section 4 we
present the numerical experiments. Section 5 concludes. Proofs of the results
 are delegated to an appendix.

\section{The model}

In this paper we consider a multi-period stochastic model of investment.
The randomness is driven by a two dimensional discrete time random walk and a Markov chain.
Let  $(b^1,b^2):=(b_{t_n}^1,b_{t_n}^2)_{n=0,1,...,\infty},$ be a
two-dimensional binomial random walk on a complete
probability space $(\Omega,\mathcal {F},\{{\mathcal {F}}_{t_n}\},
\mathbb{P})$. The random walk is assumed symmetric under $\mathbb{P}$ in the sense that
\begin{equation}\label{*0}
\P(\Delta
b_{t_n}^i=\pm 1)=1/2,\qquad i=1,2.
\end{equation}
Here $\Delta$ is the difference operator, i.e.,
$$\Delta b_{t_n}:=b_{t_{n+1}}-b_{t_n}.$$ A discrete time finite state homogeneous
Markov chain (MC) $J:=(J_{t_n})_{n=0,1,...,\infty}$ is defined on  $(\Omega,\mathcal {F},\{{\mathcal {F}}_{t_n}\},
\mathbb{P})$ and it takes values in the state space $\sS=\{{\bf 0},{\bf1}\}.$ The
$n-$step transition matrix $P^{(n)}=(P_{ij}^{(n)}),$ is defined by
$$P^{(n)}_{ij}: =\P(J_{t_n}=j|J_{t_0}=i),\ \ i,j={\bf 0},{\bf1}\qquad n=
0,1,...,\infty,$$ where $P^{(0)}_{ij}=1$ when $i=j,$ otherwise
$P^{(0)}_{ij}=0.$ We assume that the distribution of $J_0$ is known,
and $$\P(J_0=i|{{\mathcal {F}}_{0}^S})=\P(J_0=i), \qquad i={\bf
0},{\bf1}.$$ The filtration $\{  {\mathcal {F}}_{t_n} \}$ is generated by
the random walk $(b^1,b^2)$ and (MC) $J.$ The two dimensional
process $S=(S^1, S^2)$ is referred to as the forward process.

There are two securities available for trading, a riskless bond and a risky stock.
The trading horizon is $[0,T],$ with $T$ a exogenous finite horizon. We take the
bond as numeraire, thus it can be assumed to offer zero interest rate. There are $N+1$ trading dates: $t=0,h,...,Nh,$
with $\{0,h,...,Nh\}$ a partition of the interval $ [0, T].$ Here $h=T/N,$ and $t_{n}=nh,$ for $n=0,1,\cdots,N.$

 The stock price process $S^1:=\{S_t^1;t=0,h\ldots ,
(N-1)h,Nh\}$, follows the difference equation:
\begin{eqnarray}
\left\{\begin{array}{ll}\Delta S_{t_n}^1=S_{t_n}^1
(\mu^1_{t_n}h+\sigma^1_{t_n}\sqrt{h} \Delta b_{t_n}^1),\ \
 n=0,1\ldots , N-1,\\
\\S_0^1=s^1>0,\label{st1}\qquad\qquad\end{array} \right.
\end{eqnarray} for some adapted drift process $\mu^1:=\{\mu_t^1;t=0,h\ldots,
(N-1)h,Nh\}$ and volatility process $\sigma^1:=\{\sigma_t^1;t=0,h\ldots ,
(N-1)h,Nh\}$ which are chosen so that the stock price remains positive.

 In this model we assume the existence of a non-traded asset $S_t^2=\{S_t^2;t=0,h\ldots ,
(N-1)h,Nh\}$, which follows the difference equation :
\begin{eqnarray}
\left\{\begin{array}{ll}
\Delta S_{t_n}^2=S_{t_n}^2
(\mu^2_{t_n}h+\sigma^2_{t_n}\sqrt{h} \rho (\Delta b_{t_n}^1+\sqrt{1-\rho^2}\Delta b_{t_n}^2)),\ \
 n=0,1\ldots , N-1,\\
\\S_0^2=s^2>0,\label{st2}\qquad\qquad\end{array}
\right.
\end{eqnarray} for some adapted drift process $\mu^2:=\{\mu_t^2;t=0,h\ldots,
(N-1)h,Nh\},$ volatility process $\sigma^2:=\{\sigma_t^2;t=0,h\ldots,
(N-1)h,Nh\},$ and the correlation coefficient $\rho,$ with $|\rho|\leq1.$

\subsection{Trading strategies}

An investor in this model starts with initial wealth $X_0=x\in \mathbb{R}$
and trades between the stock and the bond, following self-financing strategies.
Let $\alpha_{t_n}\in {\mathcal {F}}_{t_n} $ be the wealth invested in stock at time $t_n,$ $n=0,1\ldots, N-1.$
The wealth process is governed by the self-financing equation
\begin{eqnarray}
\Delta X^{\alpha}_{t_n}&:=&\alpha_{t_n}(\mu^1_{t_n} h+\sigma^1_{t_n}
\sqrt{h} \Delta b_{t_n}^1). \label{dxt}
\end{eqnarray}
The state space of the (MC) $J$ is $\sS=\{{\bf 0},{\bf1}\},$ and it is assumed to model two states of the economy:
${\bf 0}$ the good state (bull market), ${\bf 1}$ represents the bad state (bear market). The investor's risk preference is assumed
to change according to the current state of the market
$$U(x,j)=-\exp(-\gamma(j)x)$$
where $\gamma(j),\,\,\, j\in\{{\bf 0},{\bf 1}\}$ is the positive coefficient of absolute risk aversion. 
The performance of an investment strategy is measured by the above expected utility criterion applied to the final wealth associated
with an investment strategy. The maximal expected utility is given by the value function
\begin{equation}\label{val}
 v_{t_{n}}(x,j)=\sup_{\alpha \in \Pi_{t_{n}} } \E^{\P} (- e^{-\gamma(j) X^{\alpha}_{T} } | \mathcal{F}_{t_n}).
 \end{equation}
Here $n=0,1\ldots, N-1,$ $X^{\alpha}_{T}$ is given by \eqref{dxt}, $X^{\alpha}_{t_{n}}=x,$ $J_{t_{n}}=j,$ and $\Pi_{t_{n}}$ denotes
the set of admissible trading strategies defined by
\begin{equation}\label{dd}
\Pi_{t_{n}}:=\{\alpha_{n}, \alpha_{{n+1}}, \cdots, \alpha_{{N-1}}:\alpha_{k}\in{\mathcal{F}}_{t_{k}},
\ {\rm{such\ that}}\,\,\E^{\P}|X^{\alpha}_{t_k}|<\infty, k=n,n+1\ldots N-1 \}.
\end{equation}

\subsection{A special martingale measure}
 Define a probability measures $\Q$ on
  $(\Omega,\mathcal {F},\{{\mathcal {F}}_{t_n}\})$ such that
  $(b^1,b^2)$ is a two dimensional random walk under $\Q,$
  and the following equalities holds:

$$\Q(\Delta b^1_{t_n}=1|{\mathcal {F}}_{t_n})=\frac{1-\theta^1_{t_{n}}\sqrt{h}}{2},
\ \ \ \ \Q(\Delta b^1_{t_n}=-1|{\mathcal
{F}}_{t_n})=\frac{1+\theta^1_{t_{n}}\sqrt{h}}{2}, \ n=0,1,\ldots,N-1,$$
%$$\Q(\Delta b^1_{t_n}=1|{\mathcal
%{F}}_{t_n})=\frac{1-\theta^1_{t_{n}}\sqrt{h}}{2}, \ \ \ \ \Q(\Delta
%b^1_{t_n}=-1|{\mathcal {F}}_{t_n})=\frac{1+\theta^1_{t_{n}}\sqrt{h}}{2}, \
%n=0,1,\ldots,N-1,$$
$$\Q(\Delta b^2_{t_n}=\pm1|{\mathcal
{F}}_{t_n})=\frac{1}{2},\,\,\, n=0,1,\ldots,N-1,$$
$$\Q(J_{t_n}=j|J_{t_0}=i)=\P(J_{t_n}=j|J_{t_0}=i)={P}^{(n)}_{ij}
%=\P(J_{t_n}=j|J_{t_0}=i)
\ i,j\in \sS,\,\,n=0,1,\ldots,N-1.$$  Here $\theta^1_{t_{n}},$ defined by
\begin{equation}\label{1*}
\theta^1_{t_{n}}:=\frac{\mu^1_{t_{n}}}{\sigma^1_{t_{n}}},
\end{equation}
is the market price of risk (MPR) for the stock. We choose (MPR) $\theta^1$ such that the probability measure $\Q$ is well defined, i.e.,
$$ \frac{1}{\sqrt{h}}\geq\theta^1_{t_{n}}\geq 0.$$ Expectation under probability measure $\Q$ will be denoted by $\E^{\Q}[\cdot].$
Under $\Q$ the stock price process $S^1$ is a martingale, thus we can call $\Q$ a martingale measure.
In fact this martingale measure is related to the minimal martingale measure and the minimal entropy measure
(see \cite{MZ2} ). The following Lemma is a consequence of the fact $\P$ and $\Q$
place the same probability weights on $b^2.$
\begin{lemma}\label{L11}
For every $n=0,1,...N-1,$ the following identity holds
$$\Q(S_{t_{n+1}}^2|{\mathcal
{F}}_{t_{n+1}}^{S^1}\bigvee{\mathcal {F}}_{t_n})=
\P(S_{t_{n+1}}^2|{\mathcal {F}}_{t_{n+1}}^{S^1}\bigvee{\mathcal
{F}}_{t_n}).$$
\end{lemma}

The proof is straightforward and hence omitted.

\section{Indifference valuation algorithms}

We consider a contingent claim
written on $(S^1, S^2)$ to be priced in this model.
For simplicity it is assumed that its payoff is of the form
$F(S_T^1,S_T^2)$, and is to be paid at time $T,$ with $F$ being a measurable function. Let $\lambda\geq0$ be given. The writer of $\lambda$ shares of the claim with payoff $F(S_T^1,S_T^2),$ faces the following optimization problem at time $t_n\in[0,T];$
 given the initial forward process $S_{t_n}=s=(s^1, s^2),$
the market state $J_{t_n}=i$ and the initial endowment $X_{t_n}=x,$
the writer of the option seeks an admissible trading strategy
${\bar{\alpha}}:=({\bar{\alpha}}_{k})_{k=n,n+1,...,N-1}$ such that
\begin{eqnarray}\label{*0} {\bar{\alpha}}=\arg\sup_{\alpha\in\Pi_{t_n}}\E^{\P}
[-\exp(-\gamma(i)(X^{\alpha}_T-\lambda
F(S_T^1,S_T^2))|X_{t_n}=x,J_{t_n}=i,S_{t_n}=s].\label{v}\end{eqnarray}

This optimization problem is time inconsistent due to the fact that
the risk preference is changing according to the (MC) $J.$ Indeed
the optimal strategy computed at time $t_0$ may fail to remain optimal at a later time $t_n,$ in the sense of \eqref{*0}, if $J_{t_n}\neq J_{t_0}.$
One way out of this predicament is to consider the game theoretical approach introduced by \cite{Bjo}. In order to ease notations, let
us denote $$\E_{t_n}[\cdot]:=\E^{\P}[ \cdot|X_{t_n}=x,J_{t_n}=i,S_{t_n}=s ],\,\, \E^{\P}_{t_{n}}[\cdot|\mathcal{G}]:= \E[\cdot|\mathcal{G}
\bigvee\{X_{t_{n}}=x, S_{t_{n}}=s, J_{t_{n}}=i\}]$$
$$\E^{\Q}_{t_n}[\cdot]:=\E^{\Q}[ \cdot|X_{t_n}=x,J_{t_n}=i,S_{t_n}=s ],\,\, \E^{\Q}_{t_{n}}[\cdot|\mathcal{G}]:= \E^{\Q}[\cdot|\mathcal{G}
\bigvee\{X_{t_{n}}=x, S_{t_{n}}=s, J_{t_{n}}=i\}],$$
for every
$\mathcal{G}\subset\mathcal{F}.$

\subsection {The game theoretical approach}
Since optimal trading strategies are not defined within our context, we search for subgame perfect strategies.
This is done by backward induction.  First we consider the time period
$[(N-1)h,Nh]$ (recall that $Nh=T$). At time $(N-1)h$ the writer of $\lambda$ shares of claim solve the optimization problem
\begin{eqnarray}  (P1) \qquad V_{N-1}^{\lambda}(x,i,s)=
\sup_{\alpha\in\Pi_{t_{N-1}}}
\E_{t_{N-1}}[-\exp(-\gamma(i)(X^{\alpha}_T -\lambda F(S_T^1,S_T^2))],
\label{v1}\end{eqnarray}
and this is refer to as problem (P1). In our model $\sup$ in (P1) is attained and we denote

\begin{eqnarray} \widehat{\alpha}_{{N-1}}^\lambda:=
\arg\max_{\alpha\in \Pi_{t_{N-1}}} \E[-\exp(-\gamma(J_{t_{N-1}})
(X^{\alpha}_T-\lambda F(S_T^1,S_T^2))|\mathcal{F}_{t_{N-1}}]. \label{v^11}\end{eqnarray}
Notice that in this first step we solve a pure optimization problem.
On the time period $[(N-2)h, Nh]$ one seeks subgame perfect equilibrium strategies as follows.
Let the trading strategies $\alpha$ be of the form:
\begin{equation}
\alpha_k =\begin{cases}
\widehat{\alpha}^{\lambda}_{k}, \ \ {\rm{for}}\ k={N-1},\\\\
\alpha_{k}^{\lambda},\ \ {\rm{for}}\ k=N-2,\label{alpha2}
\end{cases}
\end{equation}
for an arbitrary ${\mathcal {F}}_{t_{N-2}}-$ adapted control
$\alpha_{N-2}^{\lambda}$ such that $(\alpha^{\lambda}_{N-2},
\hat{\alpha}_{N-1}^{\lambda}) \in \Pi_{t_{N-2}}.$ The equilibrium value function is defined by
\begin{eqnarray} (P2)\qquad V_{N-2}^{\lambda}(x,i,s):
=\sup_{\alpha \in \Pi_{t_{N-2}}}
\E_{t_{N-2}}[-\exp(-\gamma(i)(X^{\alpha}_T-\lambda F(S_T^1,S_T^2))],
\label{v123}\end{eqnarray}
and this optimization problem is refer to as problem (P2). In our model $\sup$ in (P2) is attained and we denote by

\begin{eqnarray} (\widehat{\alpha}_{{N-2}}^\lambda, \widehat{\alpha}_{{N-1}}^\lambda):=
\arg \max_{\alpha\in \Pi_{t_{N-2}}}
\E[-\exp(-\gamma(J_{t_{N-2}})(X^{\alpha}_T-\lambda
F(S_T^1,S_T^2))|\mathcal{F}_{t_{N-2}}],\label{v*11}\end{eqnarray}
the equilibrium strategy at $t_{N-2}.$ Further we proceed iteratively. On the time period
$[(N-n)h,Nh]$ one restricts to trading strategies $\alpha$ of the form:
\begin{equation}
\alpha_k =\begin{cases} \widehat{\alpha}^{\lambda}_{k},
\ \ {\rm{for}}\ k=N-(n-1),N-(n-2),\cdots, N-1,\\
\\\alpha_{k}^{\lambda}\ \ {\rm{for}} \,\,\,k=N-n,
\label{alphan}\qquad\qquad\end{cases}
\end{equation}
for an arbitrary ${\mathcal {F}}_{t_{N-n}}-$ adapted control
$\alpha_{N-n}^{\lambda}$ such that
$(\alpha_{N-n}^{\lambda},\hat{\alpha}^{\lambda}_{k})_{\{k=
N-n+1,\cdots, N-1\} }\in \Pi_{t_{N-n}}.$ The equilibrium value function is defined by
\begin{eqnarray} (Pn)\qquad V_{N-n}^{\lambda}(x,i,s):=
\sup_{\alpha \in \Pi_{t_{N-n}}}
\E_{t_{N-n}}[-\exp(-\gamma(i)(X^{\alpha}_T-\lambda F(S_T^1,S_T^2))].
\label{v123}\end{eqnarray}
The $\sup$ in (Pn) is attained and we denote by

\begin{eqnarray} (\widehat{\alpha}_{{N-n}}^\lambda,
\widehat{\alpha}_{{N-n+1}}^\lambda, \cdots,
\widehat{\alpha}_{{N-1}}^\lambda):=\arg \max_{\alpha\in \Pi_{t_{N-n}}}
 \E[-\exp(-\gamma(J_{t_{N-n}})(X^{\alpha}_T-\lambda F(S_T^1,S_T^2))|\mathcal{F}_{t_{N-n}}],\label{v*11}\end{eqnarray}
the equilibrium strategy at $t_{N-n}.$ Having clarified the choice of equilibrium trading strategies, one can turn to pricing of the claim. This is done by utility
indifference; more precisely we have below the following formal definition.

\begin{definition}\label{d1}
  The writer's {\emph{indifference price}} at time
 $t_{N-n}$ of the claim with payoff $\lambda F(S_T^1,S_T^2)$  is the following Markov process
 $$p^{\lambda}_{t_{N-n}}=G^\lambda ( t_{N-n}, X^\alpha_{ t_{N-n}},J_{ t_{N-n}}, S_{ t_{N-n}}),$$ where the function $G^\lambda (t_{N-n}, x, i,s)$ is given by the equation

\begin{equation}
V_{N-n}^\lambda(x+G^\lambda (t_{N-n}, x,
i,s),i,s)=V_{N-n}^0(x,i,s).\label{p0}
\end{equation}
\end{definition}
\noindent Thus, the indifference price is the extra wealth that
makes the option writer indifferent between writing and not writing
the claim. In order to find the indifference price at time $t_{N-n}$ one, according to equation \eqref{p0}, has to find the equilibrium
value function $V_{N-n}^{\lambda}(x,i,s).$ Consequently we have to solve the problems $(P1), (P2),\cdots, (Pn).$
Let us start with $(P1).$ Here we have the following result.

\begin{theorem}\label{t1}
  The equilibrium trading strategy for problem (P1) is given by the following equations

\begin{equation}\label{33r}
 \hat{\alpha}^{\lambda}_{N-1}= \hat{\alpha}^{0}_{N-1} + H^\lambda ( t_{N-1}, X_{ t_{N-1}},J_{ t_{N-1}}, S_{ t_{N-1}}),
 \end{equation}
with
$$ \widehat{\alpha}^{0}_{N-1}=\frac{1}{2\gamma(J_{ t_{N-1}})\sigma^1_{ t_{N-1}}\sqrt{h}}
\log\bigg(\frac{1+\theta_{ t_{N-1}}\sqrt{h}}{1-\theta_{ t_{N-1}}\sqrt{h}}\bigg), $$
and
 \begin{eqnarray}\label{0e2}
 H^\lambda(t_{N-1},x,i,s)&=&\frac{1}{h \mu^1_{ t_{N-1}} \gamma(i)} \bigg[ \E_{t_{N-1}}[\log
\E_{t_{N-1}}^{\Q}[e^{\lambda \gamma(i)F(S_T^1,S_T^2)}|
{\mathcal{F}_{t_{N}}^{S^1}}]]\\\notag &-&E_{t_{N-1}}^{\Q} [\log
\E_{t_{N-1}}^{\Q}[e^{\lambda
\gamma(i)F(S_T^1,S_T^2)}|{\mathcal{F}_{t_{N}}^{S^1}}]]\bigg].
 \end{eqnarray}
Moreover, the indifference price at time $t_{N-1},$ $p^{\lambda}_{t_{N-1}}$ is given by
\begin{equation}\label{987}
p^{\lambda}_{t_{N-1}}=G^\lambda ( t_{N-1}, X_{ t_{N-1}},J_{ t_{N-1}}, S_{ t_{N-1}}),
\end{equation}
where
\begin{equation}\label{99e11}
 G^\lambda(t_{N-1},x,i,s)= \frac{1}{\gamma(i)}\E^{\Q}_{t_{N-1}} \bigg[\log
\E_{t_{N-1}}^{\Q}[e^{\lambda \gamma(i)F(S_T^1,S_T^2)}|
{\mathcal{F}_{t_{N}}^{S^1}}]\bigg].
 \end{equation}

\end{theorem}

Proof of this Theorem is done in Appendix A.

\begin{flushright}
$\square$
\end{flushright}

\noindent Next we look at the special case when $\gamma(0)=
\gamma(1)=\gamma,$ that is, the preference of the investor doesn't change
with the market state. Here we have the following result.

\begin{corollary}\label{c1}

The indifference price $p^{\lambda}_{t_{N-1}}$ is increasing in coefficient of risk aversion $\gamma.$

\end{corollary}

Proof of this Corollary is done in Appendix B.
\begin{flushright}
$\square$
\end{flushright}

Let us move to problem (P2). Here, the presentation of the results could be omitted and one can directly state the results for (Pn). We choose not to do so in order to ease the reading of the manuscript. The following result provides the indifference price and the equilibrium trading strategy.

\begin{theorem}\label{t2}

The equilibrium trading strategy for problem (P2) is given by the following equations

$$ \hat{\alpha}^{\lambda}_{N-2}= \hat{\alpha}^{0}_{N-2} + H^\lambda ( t_{N-2}, X_{ t_{N-2}},J_{ t_{N-2}}, S_{ t_{N-2}}),$$
with
$$ \widehat{\alpha}^{0}_{N-2}=\frac{1}{2\gamma(J_{ t_{N-2}})\sigma^1_{ t_{N-2}}\sqrt{h}}
\log\bigg(\frac{1+\theta_{ t_{N-2}}\sqrt{h}}{1-\theta_{ t_{N-2}}\sqrt{h}}\bigg) , $$
and
 \begin{eqnarray}\label{1e2}
 H^\lambda(t_{N-2},x,i,s)&=&\frac{1}{h \mu^1_{ t_{N-2}} \gamma(i)} \bigg[ \E_{t_{N-2}}[\log
E_{t_{N-2}}^{\Q}[\Lambda^{\lambda} _{t_{N-1} }e^{\lambda
\gamma(i)F(S_T^1,S_T^2)} |{\mathcal{F}^{ S^{1}}_{t_{N-1}}}]]\\\notag
&-&\E^{\Q}_{t_{N-2}} [\log \E_{t_{N-2}}^{\Q}[\Lambda^{\lambda}
_{t_{N-1}} e^{\lambda \gamma(i) F(S_T^1,S_T^2)}| {\mathcal{F}^{
S^{1}}_{t_{N-1}}}]]\bigg],
 \end{eqnarray}
 where $\Lambda^{\lambda} _{t_{N-1}}$ is (later) defined in \eqref{e!2}. Moreover, the indifference price at time $t_{N-2},$ $p^{\lambda}_{t_{N-2}}$ is given by
$$p^{\lambda}_{t_{N-2}}=G^\lambda ( t_{N-2}, X_{ t_{N-2}},J_{ t_{N-2}}, S_{ t_{N-2}}),$$
with
\begin{equation}\label{e11}
 G^\lambda(t_{N-2},x,i,s)= \frac{1}{\gamma(i)}\E^{\Q}_{t_{N-2}} \bigg[\log
\E_{t_{N-2}}^{\Q}[  \Lambda^{\lambda} _{t_{N-1}} e^{\lambda
\gamma(i)F(S_T^1,S_T^2)}|{\mathcal{F}^ {S^{1}}_{t_{N-1}}}]\bigg].
 \end{equation}

\end{theorem}

This Theorem is a special case of Theorem \ref{t3}.

\begin{flushright}
$\square$
\end{flushright}

Let us consider the general case, that is the problem (Pn). Here, we need to introduce the following process $(\Lambda^{\lambda}
_{t_{N-k}})$ defined by

 \begin{equation}\label{e!2}
 \Lambda^{\lambda} _{t_{N-k}}:= \frac{e^{-\gamma(J_{t_{N-k-1}}){\displaystyle{\sum_{j=N-k}^{N-1}}}\Delta
{X}_{t_j}^{\hat{\alpha}^\lambda}}}{E_{t_{N-k-1}}\left[e^{-\gamma(J_{t_{N-k-1}}){\displaystyle{\sum_{j=N-k}^{N-1}}}\Delta
{X}_{t_j}^{\hat{\alpha}^0}}\right]},\quad k=1 ,2, \cdots, N-1,
 \end{equation}
where the equilibrium strategies $(\hat{\alpha}^{\lambda}_{k})_{\{k= N-n+1\cdots N-1\} }$ were
determined from solving the problems $(P1), (P2),\ldots, P(n-1).$ The following Theorem is the main result of our paper.

\begin{theorem}\label{t3}

The equilibrium trading strategy for problem (Pn) is given by the equations
 $$ \hat{\alpha}^{\lambda}_{N-n}= \hat{\alpha}^{0}_{N-n} + H^\lambda ( t_{N-n}, X_{ t_{N-n}},J_{ t_{N-n}}, S_{ t_{N-n}}),$$
with
$$ \widehat{\alpha}^{0}_{N-n}=\frac{1}{2\gamma(J_{ t_{N-n}})\sigma^1_{ t_{N-n}}\sqrt{h}}
\log\bigg(\frac{1+\theta_{ t_{N-n}}\sqrt{h}}{1-\theta_{ t_{N-n}}\sqrt{h}}\bigg),$$
and
\begin{eqnarray}\label{2e2}
 H^\lambda(t_{N-n},x,i,s)&=&\frac{1}{h \mu^1_{ t_{N-n}}  \gamma(i)} \bigg[ \E_{t_{N-n}}
 [\log
\E_{t_{N-n}}^{\Q}[\Lambda^{\lambda}_{t_{N-n+1}}e^{\lambda
\gamma(i)F(S_T^1,S_T^2)}|{\mathcal{F}_{t_{N-n+1}}^{S^1}}]]\\
\notag&-& \E_{t_{N-n}}^{\Q}[\log
\E_{t_{N-n}}^{\Q}[\Lambda^{\lambda} _{t_{N-n+1}} e^{\lambda
\gamma(i) F(S_T^1,S_T^2)}| {\mathcal{F}_{t_{N-n+1}}^{S^1}}]]\bigg],
 \end{eqnarray}
 where $\Lambda^{\lambda} _{t_{N-n+1}}$ is defined in \eqref{e!2}. Moreover, the indifference price at time $t_{N-n},$ $p^{\lambda}_{t_{N-n}}$ is given by
$$p^{\lambda}_{t_{N-n}}=G^\lambda ( t_{N-n}, X_{ t_{N-n}},J_{ t_{N-n}}, S_{ t_{N-n}}),$$
with
\begin{equation}\label{e117}
 G^\lambda(t_{N-n},x,i,s)= \frac{1}{\gamma(i)}\E^{\Q}_{t_{N-n}} \bigg[\log
\E_{t_{N-n}}^{\Q}[  \Lambda^{\lambda} _{t_{N-n+1}} e^{\lambda
\gamma(i)F(S_T^1,S_T^2)}|{\mathcal{F}_{t_{N-n+1}}^{S^1}}]\bigg],
 \end{equation}
\begin{flushright}
$\square$
\end{flushright}

\end{theorem}

Proof of this Theorem is done in Appendix C.

\begin{flushright}
$\square$
\end{flushright}

\subsection{The case of constant CARA}

Next we look at the special case when $\gamma(0)= \gamma(1)=\gamma.$ In this special case, the equilibrium strategies coincide
with the optimal ones. Formally, we have the following results.

\begin{corollary}\label{c7}
Let us define the following value function:
\begin{eqnarray}  \bar{V}^{\lambda}_{n}(x,i,s):=
\sup_{\alpha\in \Pi_{t_{N-n}} } \E_{ {t_{N-n}}  }[-\exp(-\gamma(X^{\alpha}_T-\lambda
F(S_T^1,S_T^2))]. \label{v123no}\end{eqnarray} The $\sup$ in
 (\ref{v123no}) is attained and let
\begin{eqnarray} (\bar{\alpha}_{{N-n}}^\lambda,
\bar{\alpha}_{{N-n+1}}^\lambda, \cdots,
\bar{\alpha}_{{N-1}}^\lambda)=\arg \max_{\alpha\in \Pi_{t_{N-n}}}
 \E_{t_{N-n}}[-\exp(-\gamma(X^{\alpha}_T-\lambda F(S_T^1,S_T^2))],\label{v*11no}\end{eqnarray}
be the optimal strategy at $t_{N-n}.$ Then for every $n=0,1,\cdots,N-1,$

$$\bar{\alpha}_n^\lambda=\widehat{\alpha}_n^\lambda.$$ %and
%$$\bar{V}_k(x,i,s)=\widehat{V}_k(x,i,s),$$
\end{corollary}

Proof of this Corollary is done in Appendix D.

\begin{flushright}
$\square$
\end{flushright}

The next result provides a simple iterative scheme to compute the indifference price.

\begin{corollary}\label{c6}

The following recurrence relationship between
$p^{\lambda}_{t_{N-n}}$ and $p^{\lambda}_{t_{N-n+1}}$ holds:

$$p^{\lambda}_{t_{N-n}}=\frac{1}{\gamma }\E^{\Q}\bigg[
\log \E[e^{\gamma p_{t_{N-n+1}}^\lambda}|{\mathcal
{F}_{t_{N-n}}}\bigvee{\mathcal
{F}_{t_{N-n+1}}^{S^1}}]\bigg|{\mathcal {F}_{t_{N-n}}}\bigg].
$$
 %$$p^{\lambda}_{t_{N-n}}=\frac{1}{\gamma }E^{\Q^i}_{t_{N-n}}\log E_{t_{N-n}}[e^{\gamma
%p_{t_{N-n+1}}^\lambda}|{\mathcal {F}_{t_{N-n+1}}^{S^1}}\bigvee
%J_{t_{N-n+1}}]. $$
 Moreover, the indifference price $p^{\lambda}_{t_{N-n}}$ is increasing in $\gamma.$

\end{corollary}

Proof of this Corollary is done in Appendix E.

\begin{flushright}
$\square$
\end{flushright}

This result was proved in \cite{MZ2} for a market model in which the randomness is given by
a two dimensional random walk.

\section{ Numerical examples}

In this section we present a numerical example of our multi-period
indifference pricing model with regime switching. 
 We consider an European call option on the non-traded asset with
a strike price $K=10$. We consider the following numerical values $N =
365,\, T= 1,\, S_0^1=10,\, S_0^2=10,\, \rho=0.2,\, \mu^2 = 0.1,\,\sigma^2 = 0.50.$ The transition matrix is
$$ P^{(1)}= 
\begin{bmatrix} 
0.6 & 0.4 \\ 
0.6 & 0.4 
\end{bmatrix}. $$ 
 
In the first step we want to understand the impact a change in risk aversion has on the indifference
price. Thus, we assume that in a bull market the risk aversion is $\gamma$ and in a bear market
there is an increase of $20\%$ so that the risk aversion is $1.2\times\gamma.$

Fig 4.1 plots the percentage change in the indifference
price with the reference point being the indifference
price of a model without regime switching (neither for risk aversion nor for market coefficients). We consider that the
initial state is either bear market or bull market (the two situations are labeled on the plot). The model parameters are $\mu^1=0.1,\sigma^1 = 0.2.$

Next we explore the effect the change in stock's mean rate of return has on the indifference
price. Thus, in a bull market the stock's mean rate of return is a constant $\mu^1$ and in a bear market
there is a decrease by $0.2$ so it is $\mu^1-0.2.$

Fig 4.2 plots the percentage change in the indifference
price with the reference point being the indifference
price of a model without regime switching (neither for risk aversion nor for market coefficients).
We consider that the initial state is either bear market or bull market (the two situations are labeled on the plot).
The model parameters are  $\gamma=0.6, \sigma^1 = 0.2 .$

Finally, we want to see the effect the change in stock's volatility has on the indifference
price. Thus, in a bull market the stock's volatility is a constant $\sigma^1$ and in a bear market
there is a increase by $0.2$ so it is $\sigma^1+0.2.$

Fig 4.3 plots the percentage change in the indifference
price with the reference point being the indifference
price of a model without regime switching (neither for risk aversion nor for market coefficients).
We consider that the initiate state is either bear market or bull market (the two situations are labeled on the plot).
The model parameters are  $\gamma=0.6, \sigma^1 = 0.1.$

\section{Conclusion}
In this paper, we consider a multi-period stochastic model with regime switching.
A Markov chain which takes two possible values ( which are thought of as the two states
of markets: bull and bear) drives the regime switching. This together with a two dimensional
random walk generates the randomness of the model. In this model there is a riskless asset (bond) and two
risky assets (a stock and a non-tradable asset). This market model is incomplete, since there are
more sources of uncertainty than tradable assets.  The risky assets have stochastic mean rate of return
and volatility which can change with the state of the market. The risk preference can also change
with the state of the market. However in each state, utility is of exponential type with constant coefficient
of relative risk aversion. Due to the change in risk preference, the optimal investment strategies are
time inconsistent in this model. Indeed,  the optimal investment strategies computed at a given time
may fail to remain optimal at later times if the risk preferences are changing. In order to overcome
this predicament the subgame perfect equilibrium strategies are introduced. They are computed
iteratively by an algorithm. The goal of this paper is to price within this model a contingent claim 
written on both the stock and the non-tradable asset. Since the market is incomplete this can not be done
by a replicating portfolio argument. Therefore, we have chosen the utility indifference approach. Thus, utility
indifference prices of this contingent claim are computed via a recursive algorithm (the same algorithm
 that produces the equilibrium strategies). By running numerical experiments we wanted to understand what is
 the percentage change in price if some model parameters are changing. For example we found that the utility indifference price
 will increase by a percenatge anywhere in between  $7\%$ and $23\%$ if the market was initially in a bear (bad) state, so the coefficient of risk aversion $\gamma$
 increased by $20\%$ (this change is recorded for all values of $\gamma\in [0.1,1]$).

\section{Appendix}

\subsection{Backward Finite Difference Equations (BS$\Delta$E)}

For an arbitrary $n=0,1,\cdots,N$ let us consider the BS$\Delta$E of the form

\begin{eqnarray}
y_{t_{N-n}}=y_{t_{N-n+1}}+
f(t_{N-n},S_{t_{N-n}},J_{t_{N-n}},
{z_{t_{N-n}}})h-\sqrt{h}{z_{t_{N-n}}}{\Delta
b_{t_{N-n}}},\label{ytn1}
\end{eqnarray}
where $y_{t_{N-n+1}}$ is given and it is
${\mathcal{F}^{S}_{t_{N-n+1}}}\bigvee {\mathcal{F}^{J}_{t_{N-n}}}$
measurable. Here $z=(z^1, z^2, z^3),$ and $\Delta b=(\Delta b^1, \Delta b^2, \Delta b^3)$ with
\begin{eqnarray}
\Delta b_{t_k}^3=\left\{\begin{array}{ll}
1,\ \ \ \ {{\rm{if}}\ \Delta b_{t_k}^1\Delta b_{t_k}^2=-1},\\
-1,\ \ {{\rm{if}}\ \Delta b_{t_k}^1\Delta
b_{t_k}^2=1}.\nonumber\end{array} \right.
\end{eqnarray}
The following Lemma is a consequence of the Predictable
Representation Theorem.
\begin{lemma}\label{RepP}
The BS$\Delta$E \eqref{ytn1} has a unique $\mathcal{F}_{t_{N-n}}$ measurable solution $(y_{t_{N-n}}, z_{t_{N-n}}).$ It is given by
\begin{equation}\label{88}
z_{t_{N-n}}=\frac{1}{\sqrt{h}}\E[y_{t_{N-n+1}} \Delta b_{t_{N-n}} | \mathcal{F}_{t_{N-n}} ],
\end{equation}
$$y_{t_{N-n}}= \E[y_{t_{N-n+1}}+f(t_{N-n},S_{t_{N-n}},J_{t_{N-n}},{z_{t_{N-n}}})h | \mathcal{F}_{t_{N-n}}].$$

\end{lemma}

\subsection{Appendix A: Proof of Theorem \ref{t1} }

Let $(Y_{t_{N-1}}^{1,\lambda},{Z_{t_{N-1}}^{1,\lambda}})$ be the $\mathcal{F}_{t_{N-1}}$
solution of the following BS$\Delta$E:\\
\begin{eqnarray}
\left\{\begin{array}{ll}
Y_{t_{N-1}}^{1,\lambda}=Y_{t_N}^{1,\lambda}+f(t_{N-1},
S_{t_{N-1}},J_{t_{N-1}}, Z_{t_{N-1}}^{1,\lambda})h-\sqrt{h}
Z_{t_{N-1}}^{1,\lambda}{\Delta b_{t_{N-1}}}, \\
\\Y_{t_N}^{1,\lambda}:=\lambda F(S_T^1,S_{T}^2),\label{yt8}
\qquad\qquad\end{array} \right.
\end{eqnarray}
for some function $f(t_{N-1},\cdot)$ to be chosen later on. This solution is known to exist according to Lemma \ref{RepP}. Next we want to construct a process
$M_{t_n}^{\alpha}$ of the form
\begin{eqnarray}\label{1}
M_{t_n}^{\alpha}=-\exp(-\gamma(J_{t_{N-1}})(X_{t_n}^\alpha-Y_{t_n}^{1,\lambda})),\,\,\, n=N-1, N.
\end{eqnarray}
 Moreover it satisfies the suboptimality/optimality principle

\begin{eqnarray}\label{2}
\mbox{suboptimality}\quad
E[M_{t_N}^{\alpha}|{\mathcal{F}_{t_{N-1}}}]\leq
M_{t_{N-1}}^{\alpha}\quad \mbox{for\,\,all}\quad \alpha\in
\Pi_{t_{N-1}},
\end{eqnarray}

\begin{eqnarray}\label{3}
\mbox{optimality} \qquad
E[M_{t_N}^{\hat{\alpha}^\lambda}|{\mathcal{F}_{t_{N-1}}}]=
M_{t_{N-1}}^{\hat{\alpha}^\lambda}\quad \mbox{for\,\,some}\quad
\hat{\alpha}^\lambda\in \Pi_{t_{N-1}}.
\end{eqnarray}
Once this program is accomplished, it is easy to see that

\begin{eqnarray} \widehat{\alpha}_{{N-1}}^\lambda=
\arg\max_{\alpha\in \Pi_{t_{N-1}}} \E[-\exp(-\gamma(J_{t_{N-1}})
(X^{\alpha}_T-\lambda F(S_T^1,S_T^2))|\mathcal{F}_{t_{N-1}}]. \label{v^1111}\end{eqnarray}
Indeed, from \eqref{1}, \eqref{2} and \eqref{3} it follows that
\begin{eqnarray}\label{4}
\E_{t_{N-1}}[-\exp(-\gamma(i)(X_{t_N}^\alpha-
Y_{t_N}^{1,\lambda}))] &\leq&\E_{t_{N-1}}[-\exp(-\gamma(i)(X_{t_{N-1}}^\alpha-
Y_{t_{N-1}}^{1,\lambda}))]\\\notag&=&\E_{t_{N-1}}[-\exp(-\gamma(i)(X_{t_{N-1}}^{\widehat{\alpha}^\lambda}-
Y_{t_{N-1}}^{1,\lambda}))]\\\notag&=&\E_{t_{N-1}}[-\exp(-\gamma(i)(X_{t_{N}}^{\widehat{\alpha}^\lambda}-
Y_{t_{N}}^{1,\lambda}))]\\\label{be}&=&V_{N-1}^\lambda(x,i,s).
\end{eqnarray}
In the light of \eqref{2} and \eqref{3}
\begin{equation}\label{23}
\E_{t_{N-1}}[\exp(-\gamma(i)[(X_{t_N}^\alpha-X_{t_{N-1}}^\alpha)
-(Y_{t_N}^{1,\lambda}- Y_{t_{N-1}}^{1,\lambda})] \geq1,
\end{equation}
 with equality if
$\alpha=\widehat{\alpha}_{N-1}^\lambda$. Moreover by \eqref{dxt} and \eqref{yt8} it follows that

\begin{eqnarray}\label{poi1}\E_{t_{N-1}}\bigg[\exp(-\gamma(i)[(X_{t_N}^\alpha
-X_{t_{N-1}}^\alpha)-(Y_{t_N}^{1,\lambda}-
Y_{t_{N-1}}^{1,\lambda})])\bigg]=e^{-\gamma(i)f(t_{N-1},s,i,{z^\alpha})h}\cdot
g_{{N-1}}(\alpha,s,i,{z^\alpha}).
\end{eqnarray}
Here $z^\lambda=(Z_{t_{N-1}}^{1,\lambda}|X_{t_{N-1}}=x, S_{t_{N-1}}=s,
J_{t_{N-1}}=i),$ and
\begin{eqnarray}
g_{{N-1}}(\alpha,s,i,{z^\alpha})
&:=&E_{t_{N-1}}\bigg[e^{-\gamma(i)\alpha(\mu^1_{t_{N-1}}h+\sigma^1_{t_{N-1}}\sqrt{h}\Delta b_{t_{N-1}}^1)}\cdot e^{\gamma(i)
\sqrt{h}{z^\alpha}{\Delta
b_{t_{N-1}}}}\bigg ]\nonumber\\%&=& E_{t_{N-1}}
%\bigg[E[e^{-\gamma(i)\alpha(\mu (i)h+\sigma_1(i)\sqrt{h}\Delta
%b_{t_{N-1}}^1)}\cdot e^{\gamma(i)
%\sqrt{h}\overrightarrow{z^\lambda}\overrightarrow{\Delta
%b_{t_{N-1}}}}|{\mathcal {F}_{t_{N-1}}}\bigvee {\mathcal
%{F}_{t_{N}}^{S^1}}]\bigg]
%\nonumber\\
&=&\frac{1}{2}e^{-\gamma(i)\alpha(\mu^1_{t_{N-1}} h+\sigma^1_{t_{N-1}}\sqrt{h})}
E_{t_{N-1}} \bigg[e^{\gamma(i) \sqrt{h}{z^\alpha}{\Delta
b_{t_{N-1}}}}|{\mathcal {F}_{t_{N}}^{S^1}\bigcap A_{t_{N-1}}}\bigg]
\nonumber\\&+&\frac{1}{2}e^{-\gamma(i)\alpha(\mu^1_{t_{N-1}}h-\sigma^1_{t_{N-1}}\sqrt{h})} E_{t_{N-1}} \bigg[e^{\gamma(i)
\sqrt{h}{z^\alpha}{\Delta b_{t_{N-1}}}}|{\mathcal
{F}_{t_{N}}^{S^1}}\bigcap
A_{t_{N-1}}^c\bigg],\nonumber\end{eqnarray} with
$A_{t_{k}}:=\{ \Delta b_{t_k}^1=1\}$ and
$A_{t_{k}}^c:=\{\Delta b_{t_k}^1=-1\}.$ Direct computations show that the function $\alpha\rightarrow g_{{N-1}}(\alpha,s,i,{z^\alpha}) $ is convex, so its minimum is given by the first order condition (FOC). From the FOC we obtain the following trading strategy
\begin{eqnarray}\label{eqq}
\hat{\alpha}^{\lambda}_{N-1}&=&\frac{1}{2\gamma(J_{t_{N-1}})\sigma^1_{t_{N-1}}\sqrt{h}}
\log\bigg(\frac{1+\theta_{t_{N-1}}\sqrt{h}}{1-\theta_{t_{N-1}}\sqrt{h}}\bigg)\\\notag&+&
\frac{1}{2\gamma(J_{t_{N-1}})\sigma^1_{t_{N-1}}\sqrt{h}}\log \bigg( \frac{\E
[e^{\gamma(J_{t_{N-1}}) \sqrt{h}{ Z_{t_{N-1}}^{1,\lambda}}{\Delta b_{t_{N-1}}}}|{\mathcal
{F}_{t_{N}}^{S^1}  \bigvee\mathcal{F}_{t_{N-1}}                   }\bigcap A_{t_{N-1}}]}{\E [e^{\gamma(J_{t_{N-1}})
\sqrt{h}{Z_{t_{N-1}}^{1,\lambda}}{\Delta b_{t_{N-1}}}}|{\mathcal
{F}_{t_{N}}^{S^1}     \bigvee\mathcal{F}_{t_{N-1}}                } \bigcap A^c_{t_{N-1}}]}\bigg).
\end{eqnarray}
Next we choose $f(t_{N-1}, \cdot)$ such that
\begin{equation}\label{l11}
e^{-\gamma(J_{t_{N-1}})f(t_{N-1},S_{t_{N-1}},J_{t_{N-1}}, Z_{t_{N-1}}^{1,\lambda})h}  \,\,g_{N-1}(\hat{\alpha}^{\lambda}_{N-1},S_{t_{N-1}},J_{t_{N-1}}, Z_{t_{N-1}}^{1,\lambda})=1.
\end{equation}
Therefore, by \eqref{23} and \eqref{poi1}, it follows that $\hat{\alpha}^{\lambda}_{N-1}$ is the equilibrium
trading strategy, i.e., $\bar{\alpha}^{\lambda}_{N-1}=\hat{\alpha}^{\lambda}_{N-1}.$ Next we want to elaborate more on \eqref{eqq}. Notice that by taking conditional expectation on (\ref{yt8}), one gets
\begin{equation}\label{z}
\sqrt{h}{Z_{t_{N-1}}^{1,\lambda}}{\Delta b_{t_{N-1}}}= \lambda
F(S_T^1,S_{T}^2)- \lambda \E[F(S_T^1,S_{T}^2)|{\mathcal
{F}_{t_{N-1}}}].\end{equation} By plugging this in \eqref{eqq} it follows that
\begin{eqnarray*}\E_{t_{N-1}}[e^{\gamma(i)\sqrt{h}
{Z_{t_{N-1}}^ {1,\lambda}}{\Delta b_{t_{N-1}}}}|{\mathcal
{F}_{t_{N}}^{S^1}}\bigcap A_{t_{N-1}}]\!\!&=&\!\!
\E_{t_{N-1}}[e^{\gamma(i)\lambda F(S_T^1,S_{T}^2)}\cdot
e^{-\gamma(i)\E[\lambda F(S_T^1,S_{T}^2)|{\mathcal
{F}_{t_{N-1}}}]}|{\mathcal {F}_{t_{N}}^{S^1}}\bigcap
A_{t_{N-1}}]\nonumber\\\!\!&=&\!\! e^{-\gamma(i)\E[\lambda
F(S_T^1,S_{T}^2)|{\mathcal {F}_{t_{N-1}}}]}\cdot
\E_{t_{N-1}}[e^{\gamma(i)\lambda F(S_T^1,S_{T}^2)}|{\mathcal
{F}_{t_{N}}^{S^1}}\bigcap A_{t_{N-1}}].
\end{eqnarray*}
Consequently after some trivial manipulations, and by Lemma \ref{L11}, it follows that the equilibrium strategy can be expressed as

\begin{eqnarray}\label{eqq1}
\hat{\alpha}^{\lambda}_{N-1}&=&\frac{1}{2\gamma(J_{t_{N-1}})\sigma^1_{t_{N-1}}\sqrt{h}}
\log\bigg(\frac{1+\theta_{t_{N-1}}\sqrt{h}}{1-\theta_{t_{N-1}}\sqrt{h}}\bigg)\\\notag&+&
\frac{1}{2\gamma(J_{t_{N-1}})\sigma^1_{t_{N-1}}\sqrt{h}}\log \bigg( \frac{\E^{\Q}
[e^{\gamma(J_{t_{N-1}}) \lambda F(S_T^1,S_{T}^2) }|{\mathcal
{F}_{t_{N}}^{S^1}\bigvee\mathcal{F}_{t_{N-1}}}\bigcap A_{t_{N-1}}]}{\E^{\Q} [e^{\gamma(J_{t_{N-1}})
\lambda F(S_T^1,S_{T}^2)}|{\mathcal
{F}_{t_{N}}^{S^1}\bigvee\mathcal{F}_{t_{N-1}}}\bigcap A^c_{t_{N-1}}]}\bigg).
\end{eqnarray}
This, combined with the definition of $\Q$ imply \eqref{33r}. Let us now turn to the indifference price $p^{\lambda}_{t_{N-n}}.$ It follows from \eqref{be} that
\begin{eqnarray}\label{5.25}
V_{N-1}^\lambda(x,i,s)=e^{-\gamma{(i)}
(x-Y_{t_{N-1}}^{1,\lambda})}.
\end{eqnarray}
Moreover from the definition of the indifference price $p^{\lambda}_{t_{N-n}}$ (see (\ref{p0})) it follows that
\begin{eqnarray}\label{515}
p_{t_{N-1}}^\lambda&=&Y_{t_{N-1}}^{1, \lambda}-Y_{t_{N-1}}^{1,0}
\\\notag&=& \!\!\!\! \E[\lambda F(S_T^1,S_T^2)|{\mathcal
{F}_{t_{N-1}}}]+ f(t_{N-1},S_{t_{N-1}},J_{t_{N-1}}, Z_{t_{N-1}}^{1,\lambda})h
-f(t_{N-1},S_{t_{N-1}},J_{t_{N-1}}, Z_{t_{N-1}}^{1,0})h,
\end{eqnarray}
where the second equality comes from \eqref{yt8}. To proceed further with the computations we need the expression of $ f(t_{N-1},\cdot).$ From \eqref{l11} we get that

\begin{eqnarray}\label{el1}
f(t_{N-1},s,i,{z^\lambda})=f(t_{N-1},s,i,{z^0})+\frac{1}{\gamma(i)h}\E_{t_{N-1}}^{\Q}[\log
\E_{t_{N-1}}[e^{\gamma(i){\sqrt{h}z^\lambda}{\Delta b_{t_{N-1}}}}|{\mathcal
{F}_{t_{N}}^{S^1}}]],
\end{eqnarray}
and
\begin{eqnarray}\label{el01}
f(t_{N-1},s,i,{z^0})=\frac{1}{2}\bigg[\left(\frac{1+\theta_{t_{N-1}} \sqrt{h}}{1-\theta_{t_{N-1}}\sqrt{h}}\right)^
{-\frac{1+\theta_{t_{N-1}}\sqrt{h}}{2}}+\left(\frac{1+\theta_{t_{N-1}}\sqrt{h}}{1-\theta_{t_{N-1}}\sqrt{h}}\right)
^{\frac{1-\theta_{t_{N-1}}\sqrt{h}}{2}}\bigg].
\end{eqnarray}
Pluging \eqref{el1} and \eqref{el01} into \eqref{515} and making use of \eqref{z} we obtain \eqref{987} and \eqref{99e11}.

\begin{flushright}
$\square$
\end{flushright}

\subsection{ Appendix B: Proof of Corollary \ref{c1}}

Let us assume that $0<\gamma_1<\gamma_2;$ the following result
$$\E_{t_{N-1}}^{\Q}[e^{\lambda
\gamma_1F(S_T^1,S_T^2)}|{\mathcal {F}_{t_{N}}^{S^1}}]\leq
\E_{t_{N-1}}^{\Q}[e^{\lambda \gamma_2F(S_T^1,S_T^2)}|{\mathcal
{F}_{t_{N}}^{S^1}}]^{\frac{\gamma_1}{\gamma_2}},$$
comes from H\"older's inequality. Taking $\log$ in both sides we get
$$\frac{1}{\gamma_1}\log \E_{t_{N-1}}^{\Q}[e^{\lambda
\gamma_1F(S_T^1,S_T^2)}|{\mathcal {F}_{t_{N}}^{S^1}}]\leq
\frac{1}{\gamma_2}\log \E_{t_{N-1}}^{\Q}[e^{\lambda
\gamma_2F(S_T^1,S_T^2)}|{\mathcal {F}_{t_{N}}^{S^1}}].$$ Next, take
expectation with respect to the measure $\Q,$ and the result yields.

\begin{flushright}
$\square$
\end{flushright}

\subsection{  Appendix C: Proof of Theorem \ref{t3}}

Let  $Y_{t_{N-n+1}}^{n,\lambda}$ be defined by
\begin{eqnarray}\label{0p}e^{\gamma(J_{t_{N-n}})Y_{t_{N-n+1}}^{n,\lambda}}:=
\E\bigg[e^{-\gamma(J_{t_{N-n}}){\displaystyle{\sum_{k=N-n+1}^{N-1}}}\Delta
{X}_{t_k}^{\hat{\alpha}^\lambda}} \cdot
e^{\lambda\gamma(J_{t_{N-n}})F(S_T^1,S_T^2)}|{\mathcal
{F}^{S}_{t_{N-n+1}}}\bigvee {\mathcal
{F}^{J}_{t_{N-n}}}\bigg],
\end{eqnarray}
where $\{\hat{\alpha}^{\lambda}_{k}\}_{k=
N-n+1\cdots N-1 }$ was determined from the problems $(P1),
(P2),\ldots, P(n-1).$ Let $(Y^{n,\lambda}, Z^{n,\lambda})$ be the solution of the following BS$\Delta$E
\begin{eqnarray}
Y_{t_{N-n}}^{n,\lambda}=Y_{t_{N-n+1}}^{n,\lambda}+
f(t_{N-n},S_{t_{N-n}},J_{t_{N-n}},
Z_{t_{N-n}}^{n,\lambda})h-\sqrt{h}{Z_{t_{N-n}}^{n,\lambda}}{\Delta
b_{t_{N-n}}},\label{ytn00}
\end{eqnarray}
for some function $f(t_{N-n},\cdot)$ to be specified later. It is easy to see that
\begin{equation}\label{zn1}
\sqrt{h}{Z_{t_{N-n}}^{n,\lambda}}{\Delta b_{t_{N-n}}}=
Y_{t_{N-n+1}}^{n,\lambda}-\E[Y_{t_{N-n+1}}^{n,\lambda}|{\mathcal
{F}_{t_{N-n}}}].\end{equation}
As in (P1), we construct construct a process
$M_{t_k}^{\alpha}$ of the form
\begin{eqnarray}\label{1n}
M_{t_k}^{\alpha}=-\exp(-\gamma(J_{t_{N-n}})(X_{t_k}^\alpha-Y_{t_k}^{n,\lambda}),\,\,\, k=N-n, N-n+1.
\end{eqnarray}
 Moreover it satisfies the suboptimality/optimality principle

\begin{eqnarray}\label{2n}
\mbox{suboptimality}\quad
\E[M_{t_{N-n+1}}^{\alpha}|{\mathcal{F}_{t_{N-n}}}]\leq
M_{t_{N-n}}^{\alpha}\,\,\mbox{for\,\,all}\,\,\alpha_{N-n}^{\lambda}\,\,\mbox{such\,that}\\\notag\,\,(\alpha_{N-n}^{\lambda},\hat{\alpha}^{\lambda}_{k})_{\{k=
N-n+1,\cdots, N-1\} }\in \Pi_{t_{N-n}},   \end{eqnarray}

\begin{eqnarray}\label{3n}
\mbox{optimality} \qquad
\E[M_{t_{N-n+1}}^{\hat{\alpha}^\lambda}|{\mathcal{F}_{t_{N-n}}}]=
M_{t_{N-n}}^{\hat{\alpha}^\lambda}\,\,\mbox{for\,\,some}\,\,\hat{\alpha}_{N-n}^{\lambda}\,\,\mbox{such\,that}\\\notag\,\,(\hat{\alpha}_{N-n}^{\lambda},\hat{\alpha}^{\lambda}_{k})_{\{k=
N-n+1,\cdots, N-1\} }\in \Pi_{t_{N-n}}.
\end{eqnarray}
Then, we claim that $\hat{\alpha}_{N-n}^{\lambda}$ is an equilibrium strategy, i.e.,
$$(\widehat{\alpha}_{{N-n+1}}^\lambda, \cdots,
\widehat{\alpha}_{{N-1}}^\lambda)=\arg\max_{(\alpha_{N-n}^{\lambda},\hat{\alpha}^{\lambda}_{k})_{\{k=
N-n+1\cdots N-1\} }
\in \Pi_{t_{N-n}}}
 \E[-\exp(-\gamma(J_{t_{N-n}})(X^{\alpha}_T-\lambda F(S_T^1,S_T^2))|\mathcal{F}_{t_{N-n}}].$$
Here
$$X_{T}^\alpha:=X_{t_{N-n+1}}^{\alpha}+\sum_{k=N-n+1}^{N-1}\Delta
{X}_{t_k}^{\hat{\alpha}^\lambda},\qquad X_{T}^{\hat{\alpha}}:=X_{t_{N-n+1}}^{\hat{\alpha}}+\sum_{k=N-n+1}^{N-1}\Delta
{X}_{t_k}^{\hat{\alpha}^\lambda}.$$
Indeed, by applying iterated conditional expectation property and
 suboptimality/optimality principle we get
\begin{eqnarray*}
\E_{t_{N-n}}[-e^{-\gamma(i)(X_T^\alpha-\lambda
F(S_T^1,S_T^2))}]\nonumber
=\\\notag\E_{t_{N-n}}[-e^{-\gamma(i)X_{t_{N-n+1}}^\alpha}\cdot
\E[e^{-\gamma(i){\displaystyle{\sum_{k=N-n+1}^{N-1}}}\Delta
{X}_{t_k}^{\hat{\alpha}^\lambda}} \cdot
e^{\lambda\gamma(i)F(S_T^1,S_T^2)}|{\mathcal
{F}^{S}_{t_{N-n+1}}}\bigvee {\mathcal {F}^{J}_{t_{N-n}}}]]=\\\notag\E_{t_{N-n}}[-e^{-\gamma(i)
(X_{t_{N-n+1}}^\alpha-Y_{t_{N-n+1}}^{n,\lambda})}]
\leq\\\notag[-e^{-\gamma(i)
(X_{t_{N-n}}^\alpha-Y_{t_{N-n}}^{n,\lambda})}]
=\\\notag[-e^{-\gamma(i)( X_{t_{N-n}}^{\hat{\alpha}}-Y_{t_{N-n}}^{n,\lambda})}]
=\\\notag\E_{t_{N-n}}[-e^{-\gamma(i)(X_{t_{N-n+1}}^{\hat{\alpha}}-Y_{t_{N-n+1}}^{n,\lambda})}]
=\\\notag\E_{t_{N-n}}[-e^{-\gamma(i)X_{t_{N-n+1}}^{\hat{\alpha}}}\cdot
\E[e^{-\gamma(i){\displaystyle{\sum_{k=N-n+1}^{N-1}}}\Delta
{X}_{t_k}^{\hat{\alpha}^\lambda}} \cdot
e^{\lambda\gamma(i)F(S_T^1,S_T^2)}|{\mathcal
{F}^{S}_{t_{N-n+1}}}\bigvee {\mathcal {F}^{J}_{t_{N-n}}}]]=\\\notag
\E_{t_{N-n}}[-e^{-\gamma(i)(X_T^{\hat{\alpha}}-\lambda
F(S_T^1,S_T^2))}].
\end{eqnarray*}
Thus, $\widehat{\alpha}^\lambda_{N-n}$ is the desired equilinrium strategy on $(N-n)h$.
Recall that $$\E_{t_{N-n}}[\exp(-\gamma(i)[(X_{t_{N-n+1}}^\alpha-X_{t_{N-n}}^\alpha)-
(Y_{t_{N-n+1}}^{n,\lambda}-Y_{t_{N-n}}^{n,\lambda})]] \geq1, $$ with equality if and only if
$\alpha_{N-n}=\widehat{\alpha}^\lambda_{N-n}$. By making use of  (\ref{zn1}) and doing 
similar computations as in the case
of (P1), leads to
\begin{eqnarray}\label{eqq122}
\hat{\alpha}^{\lambda}_{N-n}&=&\frac{1}{2\gamma(J_{t_{N-n}})\sigma^1_{t_{N-n}}\sqrt{h}}
\log\bigg(\frac{1+\theta_{t_{N-n}}\sqrt{h}}{1-\theta_{t_{N-n}}\sqrt{h}}\bigg)\\\notag&+&
\frac{1}{2\gamma(J_{t_{N-n}})\sigma^1_{t_{N-n}}\sqrt{h}}\log \bigg( \frac{\E^{\Q}
[e^{\gamma(J_{t_{N-n}})   Y_{t_{N-n+1}}^{n,\lambda}                  }|{\mathcal
{F}_{t_{N-n+1}}^{S^1}\bigvee\mathcal{F}_{t_{N-n}}}\bigcap A_{t_{N-n}}]}{\E^{\Q} [e^{\gamma(J_{t_{N-n}})
   Y_{t_{N-n+1}}^{n,\lambda}}|{\mathcal
{F}_{t_{N-n+1}}^{S^1}\bigvee\mathcal{F}_{t_{N-n}}}\bigcap A^c_{t_{N-n}}]}\bigg).
\end{eqnarray}
From the definitions of $Y_{t_{N-n+1}}^{n,\lambda}$(see \eqref{0p}) and $\Q,$ after some trivial manipulations, \eqref{2e2} yields.
Following the same argument as in (P1) one gets $f(t_{N-n}, \cdot)$ to be

\begin{eqnarray}\label{0el1}
f(t_{N-n},s,i,{z^\lambda})=f(t_{N-n},s,i,{z^0})+\frac{1}{\gamma(i)h}\E_{t_{N-n}}^{\Q}[\log
\E_{t_{N-n}}[e^{\gamma(i)\sqrt{h}{z^\lambda}{\Delta b_{t_{N-n}}}}|{\mathcal
{F}_{t_{N-n+1}}^{S^1}}]],
\end{eqnarray}
and
\begin{eqnarray}\label{0el01}
f(t_{N-n},s,i,{z^0})=\frac{1}{2}\bigg[\left(\frac{1+\theta_{t_{N-n}} \sqrt{h}}{1-\theta_{t_{N-n}}\sqrt{h}}\right)^
{-\frac{1+\theta_{t_{N-n}}\sqrt{h}}{2}}+\left(\frac{1+\theta_{t_{N-n}}\sqrt{h}}{1-\theta_{t_{N-n}}\sqrt{h}}\right)
^{\frac{1-\theta_{t_{N-n}}\sqrt{h}}{2}}\bigg].
\end{eqnarray}
Furthermore, as in (P1) we get
\begin{equation}\label{tri}
p_{t_{N-n}}^\lambda=Y_{t_{N-n}}^{n,\lambda}-Y_{t_{N-n}}^{n,0}.
\end{equation} Consequently $p^{\lambda}_{t_{N-n}}=G^\lambda ( t_{N-n}, X_{ t_{N-n}},J_{ t_{N-n}}, S_{ t_{N-n}}),$ where
\begin{eqnarray}\label{g}
G^\lambda({t_{N-n}},x,i,s)&=&\E[Y_{t_{N-n+1}}^{n, \lambda}
-Y_{t_{N-n+1}}^{n, 0}|{\mathcal {F}_{t_{N-n}}}]+ f(t_{N-n},x,s,i,{z^
\lambda})h-f(t_{N-n},x,s,i,{z^
0})h\nonumber\\
&=&\frac{1}{\gamma(i)}\E_{t_{N-n}}^{\Q}\log(\frac{\E_{t_{N-n}}[e^{\gamma(i)
Y_{t_{N-n+1}}^{n,\lambda}}|{\mathcal
{F}_{t_{N-n+1}}^{S^1}}]}{\E_{t_{N-n}}[e^{\gamma(i)Y_{t_{N-n+1}}^{n,
0}} |{\mathcal
{F}_{t_{N-n+1}}^{S^1}}]})\\
&=&\frac{1}{\gamma(i)}\E_{t_{N-n}}^{\Q}\log
\frac{\E_{t_{N-n}}\bigg[e^{-\gamma(i)\displaystyle{\sum_{k=N-n+1}^{N-1}}\Delta
{X}_{t_{k}}^{\hat{\alpha}^\lambda}}e^{\lambda
\gamma(i)F(S_T^1,S_T^2)}\bigg|{\mathcal
{F}_{t_{N-n+1}}^{S^1}}\bigg]}{\E_{t_{N-n}}
\bigg[e^{-\gamma(i)\displaystyle{\sum_{k=N-n+1}^{N-1}}\Delta
{X}_{t_{k}}^{\hat{\alpha}^0}}\bigg]}\nonumber\\
&=&\frac{1}{\gamma(i)}\E_{t_{N-n}}^{\Q}\log
\E_{t_{N-n}}\left[\displaystyle{\frac{e^{-\gamma(i)\displaystyle
{\sum_{k=N-n+1}^{N-1}}\Delta
{X}_{t_{k}}^{\hat{\alpha}^\lambda}}}{\E_{t_{N-n}}
\bigg[e^{-\gamma(i)\displaystyle{\sum_{k=N-n+1}^{N-1}}\Delta
{X}_{t_{k}}^{\hat{\alpha}^0}}\bigg]}}\cdot e^{\lambda
\gamma(i)F(S_T^1,S_T^2)}\bigg|{\mathcal
{F}_{t_{N-n+1}}^{S^1}}\right]\nonumber\\
&=&\frac{1}{\gamma(i)}E_{t_{N-n}}^{\Q}\bigg[\log
\E_{t_{N-n}}[\Lambda^\lambda_{t_{N-n+1}}e^{\lambda
\gamma(i)F(S_T^1,S_T^2)}|{\mathcal {F}_{t_{N-n+1}}^{S^1}}]\bigg]
\nonumber\\
&=&\frac{1}{\gamma(i)}\E_{t_{N-n}}^{\Q}\bigg[\log
\E_{t_{N-n}}^{\Q}[\Lambda^\lambda_{t_{N-n+1}}e^{\lambda
\gamma(i)F(S_T^1,S_T^2)}|{\mathcal {F}_{t_{N-n+1}}^{S^1}}]\bigg],
\end{eqnarray}
Here, the first equality follows from \eqref{ytn00} and \eqref{0el1}.
The second equality follows from (\ref{zn1}). The third equality follows from
  $\displaystyle{{\sum_{k=N-n+1}^{N-1}}}\Delta
{X}_{t_{k}}^{\hat{\alpha}^0}$ being independent of ${\mathcal
{F}_{t_{N-n+1}}^{S^1}},$ and from the definition of $Y_{t_{N-n+1}}^{n,\lambda}$(see \eqref{0p}).
Therefore \eqref{e117} holds.
\begin{flushright}
$\square$
\end{flushright}

\subsection{Appendix D: Proof of Corollary \ref{c7}}
We use mathematical induction to prove the corollary. It is obvious that
$$\widehat{\alpha}_{N-1}^\lambda=\bar{\alpha}_{N-1}^\lambda.$$
The induction hypothesis is that
$$\widehat{\alpha}_{k}^\lambda=\bar{\alpha}_{k}^\lambda,\quad k=(N-1), (N-2),\cdots, (N-n+1).$$
 Next, we have to show that 
 \begin{equation}\label{sh1}
 \bar{\alpha}_{N-n}^\lambda=\hat{\alpha}_{N-n}^\lambda.
 \end{equation}

According to  Dynamic Programming Principle
\begin{eqnarray} (\bar{Pn})\qquad \bar{V}_{N-n}^{\lambda}(x,i,s):=
\sup_{\alpha \in \Pi_{t_{N-n}}}
\E_{t_{N-n}}[\bar{V}_{N-n+1}(X_{N-n+1},J_{N-n+1},S_{N-n+1})].
\label{v123os}\end{eqnarray}

From \eqref{v123no} and induction hypothesis it follows that
\begin{eqnarray*}
\E_{t_{N-n}}[\bar{V}_{N-n+1}(X_{N-n+1},J_{N-n+1},S_{N-n+1})]\nonumber
=\E_{t_{N-n}}[-e^{-\gamma X_{t_{N-n+1}}^\alpha}\cdot
e^{-\gamma{\displaystyle{\sum_{k=N-n+1}^{N-1}}}\Delta
{X}_{t_k}^{\bar{\alpha}^\lambda}} \cdot
e^{\lambda\gamma  F(S_T^1,S_T^2)}]\nonumber\\
=\E_{t_{N-n}}[-e^{-\gamma X_{t_{N-n}}^\alpha}\cdot
e^{-\gamma\alpha(\mu^1_{t_{N-n}} h+\sigma^1_{t_{N-n}} \sqrt{h}
\Delta b_{t_{N-n}}^1)}\cdot
e^{-\gamma{\displaystyle{\sum_{k=N-n+1}^{N-1}}}\Delta
{X}_{t_k}^{\bar{\alpha}^\lambda}} \cdot
e^{\lambda\gamma F(S_T^1,S_T^2)}]\nonumber\\
=\E_{t_{N-n}}[-e^{-\gamma X_{t_{N-n}}^\alpha}\cdot
e^{-\gamma\alpha(\mu^1_{t_{N-n}} h+\sigma^1_{t_{N-n}} \sqrt{h}
\Delta b_{t_{N-n}}^1)}\cdot
e^{-\gamma{\displaystyle{\sum_{k=N-n+1}^{N-1}}}\Delta
{X}_{t_k}^{\hat{\alpha}^\lambda}} \cdot
e^{\lambda\gamma F(S_T^1,S_T^2)}].
\end{eqnarray*}
In the light of this and \eqref{v123os}, the identitity \eqref{sh1} follows.

\begin{flushright}
$\square$
\end{flushright}

\subsection{Appendix E: Proof of Corollary \ref{c6}}
From problem (P1), we know that
\begin{equation}\label{4.61}\E[e^{-\gamma \Delta{X}_{t_{N-1}}^{\hat{\alpha}^\lambda}}\cdot
e^{\gamma\lambda F(S_T^1,S_T^2)}|{\mathcal {F}_{t_{N-1}}}]=e^{\gamma
Y_{t_{N-1}}^{1,\lambda}}.\end{equation}
From the definition of $Y_{t_{N-n+1}}^{n,\lambda}$(  with $n=2$),
 equation \eqref{4.61}, and from the iterated conditioning property we get
\begin{eqnarray}e^{\gamma Y_{t_{N-1}}^{2,\lambda}}:&=&
\E[e^{-\gamma\Delta{X}_{t_{N-1}}^{\hat{\alpha}^\lambda}}\cdot
e^{\gamma\lambda F(S_T^1,S_T^2)}|{\mathcal {F}_{t_{N-1}}^{S}}\bigvee
{\mathcal
{F}^{J}_{t_{N-2}}}]\nonumber\\
&=&\E\bigg[
\E[e^{-\gamma\Delta{X}_{t_{N-1}}^{\hat{\alpha}^\lambda}}\cdot
e^{\gamma\lambda F(S_T^1,S_T^2)}|{\mathcal
{F}_{t_{N-1}}}]\bigg|{\mathcal {F}_{t_{N-1}}^{S}}\bigvee
{\mathcal {F}^{J}_{t_{N-2}}}\bigg]\nonumber\\
&=&\E[e^{\gamma Y_{t_{N-1}}^{1,\lambda}}|{\mathcal
{F}_{t_{N-1}}^{S}}\bigvee {\mathcal {F}^{J}_{t_{N-2}}}].\nonumber
\end{eqnarray}
From the optimality principle and from the definition of $Y_{t_{N-n+2}}^{n-1,\lambda}$ if follows that
\begin{eqnarray}e^{\gamma Y_{t_{N-n+1}}^{n-1,\lambda}}&=&
\E[e^{-\gamma\Delta{X}_{t_{N-n+1}}^{\hat{\alpha}^\lambda}}\cdot
e^{\gamma Y_{t_{N-n+2}}^{n-1,\lambda}}|{\mathcal {F}_{t_{N-n+1}}}]
\nonumber\\&=&
\E[e^{-\gamma\displaystyle{\sum_{k=N-n+1}^{N-1}}\Delta{X}_{t_{k}}^{\hat{\alpha}^\lambda}}\cdot
e^{\gamma\lambda F(S_T^1,S_T^2)}|{\mathcal
{F}_{t_{N-n+1}}}].\label{mq1}
\end{eqnarray}
From the definition of $Y_{t_{N-n+1}}^{n,\lambda},$ equation
\eqref{mq1}, and from the iterated conditioning property we get
\begin{eqnarray}e^{\gamma Y_{t_{N-n+1}}^{n,\lambda}}&=&
\E[e^{\gamma Y_{t_{N-n+1}}^{n-1,\lambda}}|{\mathcal
{F}_{t_{N-n+1}}^{S}}\bigvee {\mathcal
{F}^{J}_{t_{N-n}}}],\quad n=1,2\cdots N.\label{5.32}
\end{eqnarray}
Thus, by conditioning and making use of iterated conditioning property we get
\begin{eqnarray}\label{5.36}
\E[e^{\gamma
Y_{t_{N-n+1}}^{n,\lambda}}|{\mathcal {F}_{t_{N-n+1}}^{S^1}}\bigvee
{\mathcal {F}_{t_{N-n}}}]
&=&\E\bigg[\E[e^{\gamma Y_{t_{N-n+1}}^{n-1,\lambda}}|{\mathcal
{F}_{t_{N-n+1}}^{S}}\bigvee {\mathcal
{F}^{J}_{t_{N-n}}}]\bigg|{\mathcal {F}_{t_{N-n+1}}^{S^1}}\bigvee
{\mathcal {F}_{t_{N-n}}}\bigg]
\nonumber\\
&=&\E[e^{\gamma Y_{t_{N-n+1}}^{n-1,\lambda}}|{\mathcal
{F}_{t_{N-n+1}}^{S^1}}\bigvee {\mathcal {F}_{t_{N-n}}}].%\nonumber%\\
%&=&E_{t_{N-n}}[e^{\gamma Y_{t_{N-n+1}}^{n-1,\lambda}}|{\mathcal
%{F}_{t_{N-n+1}}^{S^1}}]\nonumber
\end{eqnarray}
Therefore
\begin{eqnarray*}
p^{\lambda}_{t_{N-n}}&=&\frac{1}{\gamma}\E^{\Q}
[\log\frac{\E[e^{\gamma
Y_{t_{N-n+1}}^{{n-1},\lambda}}|{\mathcal
{F}_{t_{N-n+1}}^{S^1}} \bigvee {\mathcal {F}_{t_{N-n}}}    ] }{\E[e^{\gamma
Y_{t_{N-n+1}}^{n-1,0}}| {\mathcal {F}_{t_{N-n+1}}^{S^1}} \bigvee
{\mathcal {F}_{t_{N-n}}}]} |{\mathcal {F}_{t_{N-n}}} ] \nonumber\\
%&=&\frac{1}{\gamma}E_{t_{N-n}}^{\Q^i}\log({E_{t_{N-n}}[e^{\gamma
%Y_{t_{N-n+1}}^{n-1,\lambda}}|{\mathcal {F}_{t_{N-n+1}}^{S^1}}\bigvee
%J_{t_{N-n+1}}]})- \frac{1}{\gamma}E_{t_{N-n}}^{\Q^i}\log
%e^{\gamma(i)
%Y_{t_{N-(n-1)}}^{n-1,0}}\nonumber\\
&=&\frac{1}{\gamma}\E^{\Q}[\log\bigg({\E[e^{\gamma
(Y_{t_{N-n+1}}^{n-1,\lambda}-Y_{t_{N-n+1}}^{n-1, 0})}|{\mathcal
{F}_{t_{N-n+1}}^{S^1}} \bigvee
{\mathcal {F}_{t_{N-n}}}]}\bigg)| {\mathcal {F}_{t_{N-n}}}]\nonumber\\
&=&\frac{1}{\gamma }\E^{\Q}[\log \E[e^{\gamma
p_{t_{N-n+1}}^\lambda}|{\mathcal {F}_{t_{N-n+1}}^{S^1}}\bigvee
{\mathcal {F}_{t_{N-n}}}]|{\mathcal {F}_{t_{N-n}}}].\nonumber
\end{eqnarray*}
Here the first equality comes from \eqref{g} and \eqref{5.36}.
The second equality holds due to $Y_{t_{N-n+1}}^{n-1, 0}\in
{\mathcal {F}_{t_{N-n+1}}^{S^1}}\bigvee{\mathcal{F}_{t_{N-n}}}.$
The third equality follows from \eqref{tri}. The monotonocity of $p^{\lambda}_{t_{N-n}}$ in $\gamma$ follows from H\"older's inequality
and an induction argument as in the case of Corollary's \ref{c1} proof.

\begin{flushright}
$\square$
\end{flushright}

\end{document}